\documentclass[11pt,reqno]{amsart}
\usepackage{amsmath,amssymb,amsthm,mathrsfs}
\usepackage[top=1in,bottom=1in,left=1in,right=1in]{geometry}
\usepackage[expansion=false]{microtype}
\usepackage[dvipsnames]{xcolor}
\usepackage{xurl}
\usepackage[colorlinks=true,linkcolor=NavyBlue,citecolor=OliveGreen,urlcolor=NavyBlue]{hyperref}
\newtheorem{maintheorem}{Theorem}

\newtheorem{theorem}{Theorem}[section]
\newtheorem{proposition}[theorem]{Proposition}
\newtheorem{lemma}[theorem]{Lemma}
\newtheorem{corollary}[theorem]{Corollary}
\theoremstyle{definition}
\newtheorem{definition}[theorem]{Definition}
\theoremstyle{remark}
\newtheorem{remark}[theorem]{Remark}
\newtheorem*{remark*}{Remark}
\numberwithin{equation}{section}
\newcommand{\Q}{\mathbb Q}
\newcommand{\Z}{\mathbb Z}
\newcommand{\R}{\mathbb R}
\newcommand{\F}{\mathbb F}
\newcommand{\p}{\mathfrak p}
\newcommand{\cO}{\mathcal O}
\newcommand{\hh}{\widehat h}
\newcommand{\leg}[2]{\bigl(\tfrac{#1}{#2}\bigr)}
\DeclareMathOperator{\Nm}{N}
\DeclareMathOperator{\Num}{Num}
\DeclareMathOperator{\sqf}{sqf}
\DeclareMathOperator{\Sel}{Sel}
\DeclareMathOperator{\Cl}{Cl}
\DeclareFontFamily{U}{wncy}{}
\DeclareFontShape{U}{wncy}{m}{n}{<->wncyr10}{}
\DeclareSymbolFont{mcy}{U}{wncy}{m}{n}
\DeclareMathSymbol{\Sha}{\mathord}{mcy}{"58}
\hypersetup{pdftitle={A Diophantine definition of the integers in the rationals},pdfauthor={Manjul Bhargava and Ashay Burungale}}
\title{A Diophantine definition of $\Z$ in $\Q$}
\author{Manjul Bhargava}
\address{Department of Mathematics, Princeton University, Princeton, NJ 08544, USA}
\email{bhargava@math.princeton.edu}
\author{Ashay Burungale}
\address{Department of Mathematics, The University of Texas at Austin, Austin, TX 78712, USA}
\email{ashayk@utexas.edu}
\dedicatory{To Di Yin}
\date{November 2026}
\begin{document}
\begin{abstract}
We prove that $\Z$ is a Diophantine subset of $\Q$. More precisely, there exist an integer $m\ge1$ and a polynomial of total degree at most four in $m+1$ variables whose rational zero set projects onto $\Z$ under the first coordinate.

We construct a Diophantine model of $\Z$ on an elliptic curve and use it to define real intervals with rational endpoints, each containing an integer. A denominator bound identifies a rational number in a sufficiently short interval with that integer. The key ingredients are a criterion for square denominators, obtained using the Tian--Yuan--Zhang genus-point criterion for congruent number elliptic curves, and a height bound from quaternionic abelian surfaces. The Green--Tao theorem supplies prime progressions needed to define squaring, and hence multiplication, in the model. 
Our proof is inspired by OpenAI's approach to Hilbert's tenth problem over $\Q$, but does not rely on the $2$-converse or $2$-adic Fontaine--Mazur results announced and invoked by OpenAI.
\end{abstract}
\maketitle
\tableofcontents

\section{Introduction}

Hilbert's tenth problem asks for an algorithm deciding whether a polynomial equation with integer coefficients has an integer solution. For a commutative ring $R$, a subset of $R^r$ is \emph{Diophantine over $R$} if it is the projection of the set of $R$-solutions of a finite system of polynomial equations over $R$. The theorem of Davis, Putnam, Robinson and Matiyasevich (DPRM) identifies the Diophantine subsets of $\Z^r$ with the recursively enumerable sets\footnote{A set is recursively enumerable if an algorithm lists exactly its elements, possibly with repetitions. It is recursive, or decidable, if membership can be decided by an algorithm that always terminates.} \cite{DPR61,Mat70}. Such sets need not be decidable, so Hilbert's tenth problem over $\Z$ has a negative answer. For rational solutions, the corresponding question is Hilbert's tenth problem over $\Q$. We prove that $\Z$ is Diophantine over $\Q$ (see Theorem~\ref{thm:Z}). Replacing each integer variable by rational variables satisfying the defining equations for $\Z$ reduces integer solvability to rational solvability.

The construction represents integers by points on an elliptic curve, whose coordinates also describe the real intervals that enter the definition in Theorem~\ref{thm:Z}. Put
\[
F=\Q(\sqrt{13}),\qquad d=\frac{1+\sqrt{13}}2,\qquad
E\colon y^2=x^3-d^2x,
\]
and let $O$ be the identity of $E$. The group $\Gamma=2E(F)$ is infinite cyclic (see Lemma~\ref{lem:curve}).\footnote{This choice gives a rank-one group with the real density needed for short intervals. Complex multiplication by $\Z[i]$ and the quadratic-twist relation with the conjugate curve permit the local comparison used to define squaring.} For $P\in\Gamma\setminus\{O\}$, put
\begin{equation}\label{eq:coordinates}
u(P)=\Nm_{F/\Q}\!\left(\frac1{x(P)}\right),\qquad
y(P)=a(P)+b(P)\sqrt{13},
\end{equation}
where $u(P),a(P),b(P)\in\Q$. We have $0<u(P)<1/3$, and the pairs $(a(P),b(P))$ are dense in $\R^2$ (see Lemma~\ref{lem:density}). For $u\in\Q^\times$, let $\Num(u)$ be the absolute value of its numerator in lowest terms. We shall employ the sextic
\begin{equation}\label{eq:sextic}
G(X)=(X^2-1)(X^2-9)(X^2-25).
\end{equation}

\begin{definition}\label{def:intervals}
Fix an integer $C_0>1$. Let $T_{C_0}\subset\Q^3$ consist of the triples
\[
(u,a,b)=(u(P),a(P'),b(P')),\qquad P,P'\in\Gamma\setminus\{O\},
\]
such that the interval $[a,b]$ contains an integer and
\[
b-a<C_0^{-\Num(u)^{24}}.
\]
When $C_0$ is fixed, we write $T=T_{C_0}$.
\end{definition}

\begin{maintheorem}\label{thm:Z}
$\Z$ is a Diophantine subset of $\Q$. More precisely, the set $T_{C_0}$ is Diophantine over $\Q$ for any fixed integer $C_0>1$. If $C_0$ is sufficiently large\footnote{One may take any integer $C_0$ with $\log C_0>H+\log16$, where $H$ is the constant in Proposition~\ref{prop:height}. This choice is independent of $q$.}, a rational number $q$ is an integer if and only if, putting $g=G(16q+8)$, there exist $(u,a,b)\in T_{C_0}$ and $s,t,z\in\Q$ such that
\begin{equation}\label{eq:final-definition}
a\le q\le b,\qquad
(g-u^2)(1-g^2s^4)(1-g^2t^4)=gz^2.
\end{equation}
\end{maintheorem}

\begin{corollary}\label{cor:main}\leavevmode
\begin{enumerate}
\item A subset of $\Q^r$ is Diophantine if and only if it is recursively enumerable.
\item There exist $m\ge1$ and $\Phi\in\Z[X,W_1,\ldots,W_m]$ of total degree at most four such that $q\in\Z$ if and only if $\Phi(q,\mathbf w)=0$ for some $\mathbf w\in\Q^m$.
\item For some $N$, rational solvability of quartics with integer coefficients in $N$ variables is undecidable.
\item If $K\subseteq L$ are fields elementarily equivalent to $\Q$, then $K$ is relatively algebraically closed in $L$.
\item No algorithm decides whether a smooth projective geometrically integral variety over $\Q$ has a rational point.
\end{enumerate}
\end{corollary}
\noindent
Part~(4) answers \cite[Question~25]{Koe16}. We recall the deductions in \S\ref{sec:integers}; Remark~\ref{rem:variables} records the number of auxiliary variables once a quadratic presentation of $T$ is fixed.

\begin{remark*}
Theorem~\ref{thm:Z} also provides a counterexample to 
Mazur's conjecture that, for every
variety $V$ over $\Q$, the closure of $V(\Q)$ in $V(\R)$ has finitely
many connected components \cite{Maz92,CZ00}. Specifically, take $V$ to be the affine hypersurface $\Phi(X,\mathbf W)=0$ of Corollary~\ref{cor:main}(2).
The coordinate~$X$ maps $V(\Q)$ onto $\Z$. Since $\Z$ is closed in
$\R$, $X$ is integer-valued on the closure of $V(\Q)$ and is therefore
constant on each connected component. Every integer occurs, so there
are infinitely many such components.
\end{remark*}

\subsection{Outline of the proof}
Fix a generator $P_0$ of $\Gamma$. We represent an integer $n$ by the point $nP_0$ and call $n$ the \emph{index} of this point. The group law adds indices, and we can express multiplication through addition and squaring, since $2mn=(m+n)^2-m^2-n^2$. Our main step is therefore to show that the relation $n_2=n_1^2$ between the indices of two points is Diophantine (Proposition~\ref{prop:model}). An algorithm can decide the conditions defining $T$ from the indices of $P$ and $P'$, so we then obtain from DPRM that $T$ is Diophantine. We use two results both in the proof of Proposition~\ref{prop:model} and in the proof of Theorem~\ref{thm:Z}.

\textbf{Square denominators} (\S\ref{sec:parity}). We say that a rational number has a \emph{square denominator} if all its negative valuations are even. For $c\in\Q^\times$, we show that if the equation
\begin{equation}\label{eq:intro-parity}
(1-c)(s^4-c^2)(t^4-c^2)=z^2
\end{equation}
has a rational solution, then $c$ has a square denominator. Conversely, we show that a rational solution exists if $c$ has a square denominator, $0<c<1$, $v_2(c)\ge3$, and the squarefree part of $c$ is congruent to $7$ modulo $8$ (Lemma~\ref{lem:parity}). For the converse we need rational points with $x$-coordinate in $c\Q^{\times2}$ on two quadratic twists of the congruent number curve $y^2=x^3-x$. We construct such twists in Theorem~\ref{thm:paired}: a $2$-descent shows that points of this kind exist as soon as the rank is positive, and we deduce from the genus-point criterion of Tian, Yuan and Zhang \cite{TYZ17} that the rank is positive.

\textbf{Heights} (\S\ref{sec:height}). For $\xi\in\Q$ with $G(\xi)\ne0$, let $M(\xi)$ be the product of the primes $p>25$ at which $G(\xi)$ has positive odd valuation. We prove that the absolute logarithmic height of $\xi$ satisfies $h(\xi)\le HM(\xi)^{24}$ for a constant $H$ (Proposition~\ref{prop:height}). To do so, we attach to the point $(\xi,\sqrt{G(\xi)})$ of the curve $w^2=G(x)$ an abelian surface with quaternionic multiplication, by means of a map from this curve to a quotient of the Shimura curve of discriminant six. We bound the Faltings height of this surface by a theorem of von K\"anel--Kret \cite{vKK23}, in terms of the primes at which the surface has bad reduction or its field of definition is ramified. These primes belong to a fixed finite set or divide $M(\xi)$, and we obtain the bound for $h(\xi)$ by comparing heights.

\textbf{Squaring} (\S\ref{sec:model}). Let $P_1=n_1P_0$ and $P_2=n_2P_0$, and let $[i](x,y)=(-x,iy)$ denote complex multiplication by $i$ on $E$. For a prime $p>25$ and a place $\p$ of $F(i)$ above $p$, we give conditions on auxiliary points in $\Gamma+[i]\Gamma$ which imply that $p^{24}$ divides $n_1^2-n_2$ (Lemma~\ref{lem:local}); these conditions compare formal logarithms at $\p$ on $E$ and on its conjugate curve. If $n_2=n_1^2$ and $p\equiv17\pmod{312}$, we construct auxiliary points satisfying the conditions by means of the Frobenius endomorphism at $\p$.

For $\xi\in\Q$ with $G(\xi)\ne0$, we combine these local conditions with an instance of \eqref{eq:intro-parity} into a single Diophantine condition $\mathcal C(P_1,P_2;\xi)$ (Definition~\ref{def:globalconditions}). We show that it implies that $M(\xi)^{24}$ divides $n_1^2-n_2$ (Lemma~\ref{lem:globaldivisibility}) and hence, by the height bound, that $h(\xi)\le H|n_1^2-n_2|$ if $n_2\ne n_1^2$. We then require $\mathcal C(P_1,P_2;\xi)$ for $\xi=s$ and for $\xi=s+\beta$, where $s$ is a rational number and $\beta$ runs through the four rational coordinates of $x(m_0P_1)$ and $x(m_0P_2)$, for a fixed large integer $m_0$. If $n_2\ne n_1^2$, this bounds the heights of these coordinates by a constant multiple of $n_1^2+n_2^2$, where the constant does not depend on $m_0$; this contradicts the growth of the canonical height of $m_0P_j$ with $m_0$. Hence the five conditions imply $n_2=n_1^2$ (Lemma~\ref{lem:model-soundness}).

Conversely, if $n_2=n_1^2$, we show that the five conditions can be satisfied for a suitable $s$ (Lemma~\ref{lem:model-completeness}). We find the primes $p\equiv17\pmod{312}$ needed for this by means of the Green--Tao theorem \cite{GT08}. Since $G$ splits into linear factors over $\Q$, after clearing denominators all the shifted factors can be realized at prescribed positions in a single progression of primes.

\textbf{The integers} (\S\ref{sec:integers}). Let $q=m/n$ in lowest terms, with $n>0$, and suppose that \eqref{eq:final-definition} holds. After rescaling the variables, the equation in \eqref{eq:final-definition} becomes \eqref{eq:intro-parity} with $c=u^2/g$, so $c$ has a square denominator. At a prime dividing $M(16q+8)$, the valuation of $c$ is odd and hence positive, so the prime divides $\Num(u)$. We deduce from the height bound that $\log(n/16)\le h(16q+8)\le H\Num(u)^{24}$. For $C_0$ sufficiently large, the interval $[a,b]$, which contains both $q$ and an integer, is then shorter than $1/n$; hence this integer is $q$. Conversely, if $q$ is an integer, we use Lemma~\ref{lem:density} to choose a point $P$ for which the equation in \eqref{eq:final-definition} is solvable by the criterion above, and a point $P'$ for which $[a(P'),b(P')]$ is a sufficiently short interval containing $q$.

\subsection{Relation to other works}
Robinson gave a first-order definition of $\Z$ in $\Q$ \cite{Rob49}; Poonen obtained a universal-existential definition \cite{Poo09}, and Koenigsmann a universal one \cite{Koe16}. Rank-one elliptic curves and approximation also occur in Poonen's Diophantine models over large subrings of $\Q$ \cite{Poo03}. His earlier work uses local approximation together with global size bounds \cite{Poo02}. In the present paper, the integer represented by a point is its index, and the interval argument gives a Diophantine definition of $\Z$ itself. Koymans--Pagano \cite{KP27} and Alp\"oge--Bhargava--Ho--Shnidman \cite{ABHS26} independently prove undecidability over the rings of integers of all number fields, using rank stability in quadratic extensions for elliptic curves and abelian varieties, respectively.

OpenAI's work \cite{OAI} derives undecidability over $\Q$ from the results announced in \cite{OAI-2c,OAI-FM}, without establishing a Diophantine definition of $\Z$. The $2$-converse is used to obtain rational points representing prescribed Selmer classes, while $2$-adic Fontaine--Mazur modularity is used in the height estimate. The rational points required here are constructed in Theorem~\ref{thm:paired} using the Tian--Yuan--Zhang criterion, and the height estimate of Proposition~\ref{prop:height} is obtained from von K\"anel--Kret's theorem. The comparison of an elliptic curve with its conjugate is taken from \cite{OAI}; the local alternatives and height comparison are based on Sun's constructions \cite[Lemma~3.1, Proposition~9.1]{Sun}. Auxiliary points in $\Gamma+[i]\Gamma$ avoid a coprimality condition on the reduction orders. The rational roots of $G$ allow us to use arithmetic progressions of primes in Lemma~\ref{lem:model-completeness}.

Other approaches to a Diophantine definition have appeared in \cite{Far26,Anon26,Gla26}, independently of the present work. The constructions in \cite{Far26,Anon26} retain hypotheses from \cite{OAI}. The manuscript \cite{Gla26} proposes an unconditional proof. We have not examined its details.

\subsection*{Notation}
We write $h$ for the absolute logarithmic Weil height, $\hh$ for the canonical height on an elliptic curve, and $\sqf(c)$ for the signed squarefree integer in $c\Q^{\times2}$. Valuations satisfy $v_p(p)=1$ and $v_p(0)=+\infty$; brackets $[m]P$ denote multiplication of a point by an integer or an endomorphism. Equations over a fixed number field are expressed over $\Q$ in a fixed basis. Conjunctions are expressed by sums of squares, disjunctions by products, $a\ne0$ by $aw=1$, and $a\ge0$ by a sum of four rational squares.

\subsection*{Use of AI tools}
The second author used ChatGPT and Claude for revision and proofreading of this paper.

\subsection*{Acknowledgements}
The second author is grateful to Alex Smith for drawing his attention to the use of the $2$-converse in the OpenAI preprint \cite{OAI}, to which the present study owes its inception. The first author is grateful to the second author for the invitation to join him on this project.  We are thankful to Benedict Gross, Chandrashekhar Khare, Christopher Skinner, Ye Tian, and Don~Zagier for many inspiring discussions over the years, and also to Scott Aaronson, Levent Alp\"oge, Ben Green, Wei Ho, Peter Koymans, Rafael von K\"anel, Carlo Pagano, Ari Shnidman, and Shou-Wu Zhang for their encouraging interest and helpful comments. This work was partially supported by NSF grants DMS-1001828 and DMS-2302064.

\section{A criterion for square denominators}\label{sec:parity}

We begin by proving the criterion for square denominators (Lemma~\ref{lem:parity}). For its converse direction, we construct rational points on a pair of quadratic twists of the congruent number curve (Theorem~\ref{thm:paired}).

\begin{definition}
A rational number has a \emph{square denominator} if its positive denominator in lowest terms is a square. Equivalently, all its negative valuations are even.
\end{definition}

For $D\ne0$, write $E_D:y^2=x^3-D^2x$, with identity $O$.

\begin{lemma}\label{lem:parity}
Let $c\in\Q^\times$. If
\begin{equation}\label{eq:parity}
(1-c)(s^4-c^2)(t^4-c^2)=z^2
\end{equation}
has a solution in $\Q$, then $c$ has a square denominator. Conversely, a solution exists if $c$ has a square denominator and
\[
0<c<1,\qquad v_2(c)\ge3,\qquad \sqf(c)\equiv7\pmod8.
\]
\end{lemma}

For the converse, we need only the rational-point assertion of the following theorem.

\begin{theorem}\label{thm:paired}
Let $j,k$ be coprime positive squarefree integers with $j\equiv7\pmod8$ and $k\equiv1\pmod8$, each a square modulo any prime dividing the other. There exist infinitely many positive squarefree integers $n\equiv5\pmod8$, divisible by $j$ and prime to $k$, such that $E_n$ and $E_{nk}$ have analytic and algebraic rank one, a finite Tate--Shafarevich group of odd order, and a rational point with $x$-coordinate in $j\Q^{\times2}$ and non-zero $y$-coordinate.
\end{theorem}

For example, the rational-point assertion holds for $j=7$, $k=113$ and $n=21$:
\[
(28,98)\in E_{21}(\Q),\qquad
(23548,3595130)\in E_{2373}(\Q).
\]
Their $x$-coordinates are $7\cdot2^2$ and $7\cdot58^2$, respectively.

\begin{proof}[Proof of Lemma~\ref{lem:parity}]
If $v_p(c)$ is negative and odd, then $1-c$ has odd valuation and both quartics have finite even valuation: the valuations of $s^4$ and $c^2$ differ, and likewise for $t$. Their product cannot be a square.

For the converse put $j=\sqf(c)$ and $k=\sqf(1-c)$. If $k=1$, write $1-c=r^2$ with $r\in\Q$ and take $(s,t,z)=(0,0,c^2r)$.
Otherwise $j,k$ are positive odd integers, with $j\equiv7$ and $k\equiv1\pmod8$; the latter follows from $1-c\in1+8\Z_2$. At $p\mid j$, one has $v_p(c)>0$, so $1-c$ is a local square. At $p\mid k$, a negative valuation of $c$ would make $v_p(1-c)$ even, whereas a positive valuation would make it zero. Thus $c$ is a unit congruent to $1$ modulo $p$. These observations give coprimality and the mutual square conditions for $j,k$.

Theorem~\ref{thm:paired} therefore provides, for $D=n,nk$, a point $(x,y)\in E_D(\Q)$ with $x=cr^2$ and $y\ne0$. Set $s_D=y/(r(x-D))$. Then
\[
s_D^2=c\frac{x+D}{x-D},\qquad
s_D^4-c^2=\frac{4c^2xD}{(x-D)^2}\in cD\Q^{\times2}.
\]
The product for $D=n,nk$ has squareclass $k$, the squareclass of $1-c$, proving \eqref{eq:parity}.
\end{proof}

We recall the part of $2$-descent needed below (see Monsky's appendix to \cite{HB94}). Let $D>0$ be squarefree with $D\equiv5\pmod8$. Put $I=\{p:p\mid D\}$. For an odd prime $p\nmid a$, let $[a/p]\in\F_2$ be $0$ or $1$ according as $a$ is a square or a nonsquare modulo $p$. Let
\[
A_{pq}=[q/p]\quad(p\ne q),\qquad A\mathbf1=0,
\qquad y_p=[-1/p],\quad z_p=[2/p].
\]
Here $\mathbf1$ is the all-ones vector, so $A\mathbf1=0$ fixes the diagonal. Write $D_y,D_z$ for the diagonal matrices with diagonals $y,z$, and $d(a)=\prod_{a_p=1}p$ for $a\in\F_2^I$.

\begin{lemma}\label{lem:descent}
The group $E_D(\Q)/(2E_D(\Q)+E_D[2])$ embeds in
\begin{equation}\label{eq:point-descent}
\{(a,b):Aa=D_yb,\quad Ab=D_yb+D_za\}.
\end{equation}
A non-zero class corresponding to $(a,b)$ has a representative with $x$-coordinate in $d(a)\Q^{\times2}$. The same injection holds for $\Sel_2(E_D/\Q)/\delta(E_D[2])$, where $\delta$ is the Kummer map.
\end{lemma}

\begin{proof}
The Kummer map $(x,x-D)$ modulo squares has kernel $2E_D(\Q)$ and takes $(0,0)$ and $(D,0)$ to $(-1,-D)$ and $(D,2)$, respectively \cite[Ch.~X, Proposition~1.4]{Sil09}. Translation by these points independently changes the sign of the first entry and the valuation parity of the second at $2$. Thus each class modulo $E_D[2]$ has a unique pair with positive first entry and even second valuation at~$2$.

A non-zero normalized class comes from a point with $x>D$. The usual support argument for $x(x-D)(x+D)=y^2$ gives even valuations for both entries outside $2D$. At $2$, the second entry has even valuation by normalization; so does the first, since $x\pm D$ are units if $v_2(x)>0$, and all three factors share the valuation of $x$ if $v_2(x)<0$. Thus both entries are represented by positive divisors $d(a),d(b)$ of $D$.

At $p\mid D$, translate by $a_p(D,0)+b_p(0,0)$. The two valuation parities become even. The first coordinate $X$ of the translated point has $v_p(X)\le0$: a positive even valuation would give $v_p(X-D)=1$. Hence $1\pm D/X$ are squares in $\Q_p$, so the three factors of $X(X-D)(X+D)$ have the same squareclass and are squares. The resulting first two classes are
\[
d(a)D^{a_p}(-1)^{b_p},\qquad d(b)D^{b_p}(-1)^{b_p}2^{a_p}.
\]
Their unit parts give \eqref{eq:point-descent}.

The same argument applies to a class in $\Sel_2(E_D/\Q)$. We normalize its pair of squareclasses as above, and we use local points representing the class over $\R$ and over each $\Q_p$ in place of the rational point. The point over $\R$ shows that the second entry is positive, and the points over $\Q_p$ give the support restrictions and the symbol equations as before. A local point may lie in $E_D[2]$, which causes no difficulty. The pairs of the two-torsion points, given above, satisfy the support restrictions. At $p\mid D$, the translated point cannot be a non-zero two-torsion point, because both valuation parities are even; if it is $O$, its pair is trivial, which again gives the symbol equations.
\end{proof}

\begin{proof}[Proof of Theorem~\ref{thm:paired}]
We construct $n$ so that both descent systems reduce to the same two-by-two system, whose non-zero solution has $a$-coordinates supported exactly on the primes of $j$.

\emph{The construction.} Let $J,K$ be the prime supports of $j,k$. Draw a star on $J$ with any centre $p_\ast\in J$, and a second star with centre $p_0\equiv3\pmod8$ and leaves $K$. The two labels at a prime vertex $p$ are its coordinates $a_p,b_p$ in \eqref{eq:point-descent}. Subdivide each edge by a prime $\ell\equiv1\pmod8$, prescribing $[q/\ell]=1$ at its two neighbours and $0$ at all other vertices $q\ne\ell$, including the other inserted primes. For $p\in J$, prescribe
\[
[p_0/p]=\deg(p)+y_p,
\]
where degrees are taken before subdivision. Choose the inserted primes successively, fixing their symbols with $J\cup K$ and the earlier inserted primes. Then choose $p_0\equiv3\pmod8$ with the displayed symbols at $J$ and with $[p_0/\ell]=0$ or $1$ according as $\ell$ belongs to the first or second star. Quadratic reciprocity, the Chinese remainder theorem and Dirichlet's theorem permit these choices, with all primes distinct and outside $J\cup K$. Put $n=jp_0\prod\ell$.

\emph{The descent systems.} Consider \eqref{eq:point-descent} for $D=n$ and $nk$. At an inserted vertex $\ell$ with neighbours $p,q$, its equations read $a_p+a_q=b_p+b_q=0$, since $y_\ell=z_\ell=0$. The original labels are therefore constant on each star. Write $(a_i,b_i)$ for the common labels on the $i$-th star. For $D=n$, a leaf $q\in K$ is absent, and its inserted row instead gives $a_\ell=a_{p_0}$ and $b_\ell=b_{p_0}$. The unrestricted symbols $[p_0/q]$, $q\in K$, multiply differences of labels within the second star and thus do not affect the reduced system.

The equations at the original leaves determine the inserted labels.
When we sum over the original vertices of either star, each label at an
inserted vertex with both neighbours present occurs twice and cancels.
For $D=n$, an inserted vertex whose leaf is absent has the same labels
as $p_0$, so its contribution to the sum also vanishes. Once the leaf
equations hold, the resulting sums are equivalent to the centre
equations and give
\begin{equation}\label{eq:two-stars}
\overline A=\begin{pmatrix}1&1\\0&0\end{pmatrix},
\qquad \overline y=(1,1),\qquad \overline z=(0,1).
\end{equation}
Indeed, the upper right entry is $\sum_{p\in J}[p_0/p]=1$, since the degree sum is even and $j\equiv3\pmod4$. Reciprocity gives the reverse entry $0$. The mutual square conditions give $[k/p]=0$ for $p\in J$ and $[j/q]=0$ for $q\in K$, so adding $K$ changes neither entry. The row sums determine the diagonal, and the residues modulo $8$ determine $\overline y,\overline z$.

The reduced equations give $b_1=a_1+a_2$, $b_2=0$, then $a_2=0$. Their solution space is a line, generated by $a_1=b_1=1$, $a_2=b_2=0$. All reduced solutions lift by the leaf rules above, so the full solution space is also a line. For its non-zero vector, the leaf equations give $a_\ell=0$ at every inserted vertex, so $a$ is the characteristic vector of $J$ and $d(a)=j$.

\emph{The genus-point criterion.} It remains to show that this class is represented by a rational point. Write $\Cl(L)$ for the ideal class group of a number field $L$ and $2\Cl(L)$ for its subgroup of squares. For a positive squarefree integer $h_0$, put
\[
\gamma(h_0)=\#\bigl(2\Cl(\Q(\sqrt{-h_0}))\bigr)\pmod2.
\]
We use the following consequence of the Tian--Yuan--Zhang genus-point criterion: for $D\equiv5\pmod8$,
\begin{equation}\label{eq:paired-genus-sum}
\sum_{\substack{D=d_0\cdots d_r\\
d_0\equiv7\ (8),\ d_1\equiv3\ (8)\\
d_i\equiv1\ (8)\ (i\ge2)}}\prod_{i=0}^r\gamma(d_i)=1
\quad\Longrightarrow\quad
E_D(\Q)\ \text{has positive rank and}\ L'(E_D,1)\ne0.
\end{equation}
Here $r\ge1$, all factors exceed $1$, and the factors indexed by $i\ge2$ are unordered. The criterion is stated with this sum by Smith \cite[Theorem~2.1 and Table~2, row~5(b)]{Smi16}. We use the second of the two sums appearing in the proof of \cite[Theorem~3.5(2), p.~760]{TYZ17}, and record the implication here. In the notation of \cite{TYZ17}, with $K_D=\Q(\sqrt{-D})$, that proof forces each of its two bracketed sums, indexed by $d_0\equiv5$ and by $d_0\equiv7\pmod8$, to be even if $P(D)\in A(K_D)^-+A[4]$. The second is the sum in \eqref{eq:paired-genus-sum}, since $D\equiv5$ and $d_0\equiv7$ force $d_1\equiv3\pmod8$. Oddness therefore puts $P(D)$ outside that group. By Lemma~3.18, a torsion $P(D)$ would lie in $A[4]$, so $P(D)$ is non-torsion.
By Theorem~3.5(1), the image of $2P(D)$ under a $2$-isogeny from $A$ to $E_1\colon y^2=x^3-x$ lies in $E_1(K_D)^-$, the subgroup on which the non-trivial automorphism of $K_D$ acts by $-1$. Twisting identifies this subgroup with $E_D(\Q)$, which therefore has positive rank.

For the analytic-rank assertion, we use the opening statement of \cite[Theorem~3.5]{TYZ17}: it identifies the class of $P(D)$ after tensoring with $\Q$ with the vector defined on p.~739, which vanishes when $\mathscr L(D)=0$. Hence $\mathscr L(D)\ne0$. By the definition on p.~722, the order of vanishing of $L(E_D,s)$ at $s=1$ is then at most one; it is odd, since the root number is $-1$. Thus $L'(E_D,1)\ne0$, and the analytic rank is one.

\emph{The genus sum.} We compute the sum in \eqref{eq:paired-genus-sum} by deleting leaves. For a positive squarefree $h_0>1$ with $h_0\equiv1,3,7\pmod8$, let $A_{h_0}$ be defined as $A$ above, using the primes dividing $h_0$, and put $z_{h_0}=([2/p])_{p\mid h_0}$. R\'edei's formula gives, for any prime divisor $v$ of $h_0$,
\begin{equation}\label{eq:redei-minors}
\gamma(h_0)=
\begin{cases}
\det A_{h_0}^{(v,v)},&h_0\equiv3,7\pmod8,\\[2pt]
\det\bigl(A_{h_0}^{(v)}\mid z_{h_0}\bigr),&h_0\equiv1\pmod8.
\end{cases}
\end{equation}
Here $(v)$ means deleting column $v$, and $(v,v)$ means deleting row and column $v$; all determinants are over $\F_2$, with empty determinant $1$. Indeed, in the first case the R\'edei matrix has rank $t-1$ exactly when $\gamma(h_0)=1$, where $t$ is the number of prime divisors of $h_0$, and both its left and right kernels contain $\mathbf1$. In the second case, the genus relations are $A_{h_0}a+\epsilon z_{h_0}=0$, with $\epsilon$ the coordinate at the ramified prime $2$; fixing $a_v=0$ gives the second determinant. See \cite{Red34} and \cite[Theorem~3.1 and (16)]{Ste22}, or \cite[Table~1]{Smi16} for the minor formulas.

An inserted prime $\ell\equiv1\pmod8$ has zero $y$- and $z$-entries, and its only non-zero off-diagonal entries are the symmetric links to its neighbours. If neither neighbour lies in a factor $h_0$ containing $\ell$, choose $v\ne\ell$ in \eqref{eq:redei-minors}; the $\ell$-row is zero, so $\gamma(h_0)=0$. If $h_0=\ell$, the second determinant is $[z_\ell]=[0]$, with the same conclusion. If exactly one neighbour $p$ lies in $h_0$, choose $v=p$. The $\ell$-row then has only its diagonal entry $1$, and expansion along that row gives
$\gamma(h_0)=\gamma(h_0/\ell).$
Only the old diagonal at $p$ changes on deleting $\ell$, and its column has already been removed.

Order the original star edges whose endpoints occur in $D$. In a non-zero summand of \eqref{eq:paired-genus-sum}, each inserted prime lies in a factor containing at least one of its neighbours. If an edge has endpoints in different factors, move the inserted prime of the first such edge between those factors. The deletion rule preserves the product of the $\gamma$'s, and the factor residues and the locations of all original vertices remain unchanged. Thus this is a fixed-point-free involution, also when the factors of residue $1$ are unordered. When $D=n$, an inserted prime whose $K$-neighbour is absent must stay in the factor containing $p_0$.

In a term not paired by the involution, the remaining original vertices of each star lie in a single factor, and no factor consists only of inserted primes. The two stars cannot share a factor, since their product is $5\pmod8$. Hence only the factorization into the two stars remains; their products are $7$ and $3\pmod8$.

For either surviving factor, delete its centre in the first determinant of \eqref{eq:redei-minors}. Each inserted prime with both endpoints present has diagonal zero and a single remaining entry $1$ in both its row and column, at its original leaf. Expanding along that row and column deletes the pair with factor $1$, regardless of the entries between original vertices. The resulting empty determinant is $1$. For $D=n$, deleting the centre of the second star instead leaves an identity matrix on the inserted primes whose $K$-neighbours are absent. Thus both surviving factors have $\gamma=1$, and the genus sum is odd.

\emph{Conclusion.} By \eqref{eq:paired-genus-sum}, $E_D(\Q)$ has positive rank and $L'(E_D,1)\ne0$ for $D=n,nk$. Positive rank gives a non-zero class modulo $2E_D(\Q)+E_D[2]$. By Lemma~\ref{lem:descent}, it corresponds to the non-zero solution found above, with $d(a)=j$, and so has a representative with $x$-coordinate in $j\Q^{\times2}$ and non-zero $y$-coordinate. This proves the rational-point assertion. Since the analytic rank is one, Gross--Zagier and Kolyvagin give algebraic rank one and finite $\Sha(E_D/\Q)$ \cite{GZ86,Kol90}, and Lemma~\ref{lem:descent} gives $\dim_{\F_2}\Sel_2(E_D/\Q)=3$. Hence $\Sha(E_D/\Q)[2]=0$, so its order is odd. Finally, $p_0$ ranges over an infinite Dirichlet set with all other primes fixed.
\end{proof}

\section{A height estimate for the split sextic}\label{sec:height}

We next bound the height of a rational number $x$ in terms of the primes $p>25$ at which $G(x)$ has positive odd valuation, where $G$ is the sextic \eqref{eq:sextic}. For $x\in\Q$ with $G(x)\ne0$, let $M(x)$ be the product of these primes.

\begin{proposition}\label{prop:height}
There exists a constant $H\ge1$ such that
\begin{equation}\label{eq:sextic-height}
h(x)\le H M(x)^{24}\qquad(x\in\Q,\ G(x)\ne0).
\end{equation}
\end{proposition}

In the proof, we attach to the point $(x,\sqrt{G(x)})$ an abelian surface with quaternionic multiplication, bound its Faltings height, and then compare that height with $h(x)$.

Let $\cO$ be a fixed maximal order in the quaternion algebra of discriminant six. Fix a positive involution on $\cO$ induced by an element of square $-6$. It determines a principal polarization on every abelian surface $A$ equipped with an action $\iota:\cO\hookrightarrow\operatorname{End}(A)$ \cite[\S1]{Buz97}. Let $X/\Q$ be their coarse Shimura curve, and let $X^*$ be its maximal Atkin--Lehner quotient. The latter is isomorphic to $\mathbb P^1$ over $\Q$ and has three rational elliptic points of orders $4,2,6$ \cite[\S3.1]{Elk98}. Choose the coordinate $t$ over $\Q$ with these points at $0,1,\infty$, respectively; thus $t=1-t_{\mathrm{Elk}}$ in the notation of \cite[\S3.1, after (33)]{Elk98}.

\begin{lemma}\label{lem:QM-isogenies}
Let $(A,\iota)$ represent a point of $X(\overline\Q)$ whose image in $X^*$ is rational. For any $\tau\in G_\Q$, there exists an isogeny $\mu_\tau:\tau^*A\to A$ satisfying
\[
\mu_\tau\circ\tau^*(\iota(b))=\iota(b)\circ\mu_\tau
\qquad(b\in\cO).
\]
\end{lemma}
\begin{proof}
The marked moduli problem and $X\to X^*$ are defined over $\Q$, so $(\tau^*A,\tau^*\iota)$ and $(A,\iota)$ differ by an Atkin--Lehner action. This is induced by quotienting $\tau^*A$ by $\tau^*A[I]=\bigcap_{b\in I}\ker\tau^*\iota(b)$ for a two-sided ideal $I$ of $\cO$ \cite[p.~10]{Elk98}. The subgroup is $\cO$-stable, since $ba\in I$ for $a\in\cO$, $b\in I$, so the quotient map is $\cO$-linear. Its composite with the marked isomorphism to $(A,\iota)$ is the required~$\mu_\tau$.
\end{proof}

\begin{proof}[Proof of Proposition~\ref{prop:height}]
\emph{The map to $X^*$.} On the smooth projective curve $C\colon w^2=G(x)$, put $Y=(x+1)/2$. Its six branch values on the $Y$-line are $-2,-1,0,1,2,3$. Consider
\begin{equation}\label{eq:shimura-map}
t=\frac{27Y^2}{4(Y^2-1)^3},\qquad
1-t=\frac{(Y^2-4)(2Y^2+1)^2}{4(Y^2-1)^3}.
\end{equation}
On the $Y$-line, $t$ has zeros of orders $2,4$ at $0,\infty$, and poles of order $3$ at $\pm1$. The zeros of $1-t$ have order $1$ at $\pm2$ and order $2$ at the roots of $2Y^2+1$. The double cover defining $C$ therefore makes the local degrees above $0,1,\infty$ equal to $4,2,6$, respectively.

\emph{A Shimura cover.} Adding level $V_1(5)$ gives a fine moduli curve over $\Q$, smooth and projective, carrying a universal marked abelian surface \cite[Theorem~2.1, Lemma~2.2 and Corollary~2.3]{Buz97}. Pass to a finite \'etale cover $S$ over which the universal surface carries a symmetric line bundle $\mathcal L$ inducing the polarization, rigidified at the zero section, together with compatible theta data for $\mathcal L^{\otimes4}$, corresponding to $r=2$ in \cite[\S2.3]{Paz12}. The cover can be taken \'etale: the symmetric representatives form a torsor under the dual surface's $2$-torsion. After trivializing the $4$-torsion, choose normalized lifts in the theta group of $\mathcal L^{\otimes4}$. Requiring a lift of a torsion-basis element to have fourth power $1$ leaves a $\mu_4$-torsor; the remaining compatible normalizations also involve only roots of unity \cite[\S6]{Mum67}. These are finite \'etale choices in characteristic zero. We keep the entire cover $S$, which is defined over $\Q$ and need not be geometrically connected, and write $(\mathcal A,\iota,\mathcal L)$ for its family.

\emph{An \'etale cover of $C$.} The groups uniformizing the components of $S(\mathbb C)$ are torsion-free, since the fine level removes stabilizers. Thus $S\to X^*$ has ramification indices $4,2,6$ above $0,1,\infty$ and is unramified elsewhere. Locally above these three points, both maps have the form $z\mapsto z^r$ with the same index $r$, so their normalized fibre product is unramified over $C$. Elsewhere $S\to X^*$ is \'etale. Thus the normalized fibre product $Z\to C$ is finite \'etale.

\emph{Lifting rational points.} Outside a fixed finite set $S_0$ of rational primes containing $2$, the curves have proper smooth models, $Z\to C$ remains finite \'etale, and the family extends to an abelian scheme with its $\cO$-action. Put $D_x=\sqf(G(x))$. Properness extends the point $(x,\sqrt{G(x)})$ over the valuation rings of $\Q(\sqrt{D_x})$ outside $S_0$. Pulling back the \'etale cover and choosing a field factor gives a lift over a field $L_x$, unramified over $\Q(\sqrt{D_x})$ outside $S_0$, with
\[
[L_x:\Q]\le2\deg(Z/C).
\]
Write $s_x\in S(L_x)$ for its image and $A=\mathcal A_{s_x}$. Set $N_0=\prod_{p\in S_0}p$ and $U_x=\operatorname{Spec}\cO_{L_x}[1/N_0]$. The surface $A$ extends to an abelian scheme over $U_x$. Let $N_{U_x}$ be the product of the norms of the omitted finite places. Writing $\operatorname{rad}(m)$ for the product of the prime divisors of a non-zero integer $m$, the ramification bound implies
\[
\operatorname{rad}\bigl(N_{U_x}|\operatorname{disc}(L_x)|\bigr)\mid N_0|D_x|.
\]

\emph{The Faltings height.} The image of $(A,\iota)$ in $X^*$ is the rational point $t(x)$. Fix an embedding $F_0=\Q(\sqrt2)\hookrightarrow\cO\otimes\Q$, as in \cite[(30)]{Elk98}. Its action is defined over $L_x$. Restricting the isogenies of Lemma~\ref{lem:QM-isogenies} to $F_0$ makes $A$ of $\mathrm{GL}_2$-type with $G_\Q$-isogenies in the sense of \cite[(7.1)]{vKK23}.

Write $h_F$ for the stable Faltings height. By \cite[Theorem~7.1(i) and the sharpening immediately following it]{vKK23}, the exponent in the height bound may be taken to be $\max\{24,5\dim A\}=24$. Since $[L_x:\Q]$ is bounded, this gives, uniformly in $x$,
\begin{equation}\label{eq:faltings-height}
h_F(A)\ll\operatorname{rad}\bigl(N_{U_x}|\operatorname{disc}(L_x)|\bigr)^{24}\ll |D_x|^{24}.
\end{equation}

\emph{Comparing heights.} We compare this bound with $h(x)$ through the theta-null map. Let $\pi:\mathcal A\to S$ be the structure morphism. The chosen theta data identify $\pi_*\mathcal L^{\otimes4}$ with a fixed $16$-dimensional vector space tensored with a line bundle on $S$ \cite[\S6, Proposition~2, pp.~80--82]{Mum67}. Projectivizing removes the latter line bundle. Evaluating the relative theta embedding at the zero section therefore gives a morphism $\Theta:S\to\mathbb P^{15}$. By \cite[Definition~2.6]{Paz12},
\[
h(\Theta(s))=h_\Theta(\mathcal A_s,\mathcal L_s)+O(1),
\qquad s\in S(\overline\Q),
\]
where the bounded term accounts for the difference between the usual projective height and the $\ell^2$-height used there.

A theta-null point determines finitely many polarized surfaces \cite[\S2.3.3]{Paz12}, each with countably many $\cO$-actions and finitely many level and theta choices. Thus every fibre of $\Theta$ in $S(\mathbb C)$ is countable, so $\Theta$ is non-constant on every geometric component of $S$, and the pullback of the hyperplane class under $\Theta$ is ample. Height comparison \cite[Theorem~B.3.2]{HS00} and \cite[Corollary~1.3(2)]{Paz12} yield
\[
h(t(x))\ll 1+h(\Theta(s_x))\ll1+\max\{0,h_F(A)\}.
\]
The normalizations of $h_F$ in the two height estimates differ by a fixed constant.

\emph{Conclusion.} By \eqref{eq:shimura-map}, $h(t(x))=6h(x)+O(1)$. Since $G$ is monic of degree six with integral coefficients, $v_p(G(x))=6v_p(x)$ whenever $v_p(x)<0$. Thus all primes $p>25$ dividing $D_x$ are counted by $M(x)$, and
\[
|D_x|\le\Bigl(\prod_{p\le25}p\Bigr)M(x).
\]
Together with \eqref{eq:faltings-height}, this proves \eqref{eq:sextic-height}.
\end{proof}

\begin{remark}\label{rem:exponent}
We use only the existence of some fixed exponent in \eqref{eq:sextic-height}. If $24$ is replaced there by an integer $e\ge24$, the arguments of \S\S\ref{sec:model}--\ref{sec:integers} remain valid after three changes: $\Num(u)^{24}$ becomes $\Num(u)^e$ in Definition~\ref{def:intervals}, the primes $p>25$ are replaced by the primes $p>e+1$, and the series \eqref{eq:truncated-log} is truncated at degree $e$. In particular, the exponent $24\dim A=48$ displayed in \cite[Theorem~7.1(i)]{vKK23} would suffice.
\end{remark}

\section{Representing integers by points on an elliptic curve}\label{sec:curve}

We now establish the group structure of $\Gamma$ and the real and $p$-adic approximation properties of its points that we need for our Diophantine definition of $\Z$.

Let $E/F$ be the curve of the introduction, and let $\sigma$ be the non-trivial automorphism of $F/\Q$.
The point $Q=(d+1,d+2)$ lies on $E$ because $d^2-d-3=0$.
We represent an affine point $(x_0+x_1\sqrt{13},y_0+y_1\sqrt{13})$ by $(1,x_0,x_1,y_0,y_1)$ and the identity $O$ by $(0,0,0,0,0)$.

\begin{lemma}\label{lem:curve}
The curve $E$ has rank one over $F$, and $E(F)_{\rm tors}=E[2]$.
Thus $\Gamma=2E(F)$ is an infinite cyclic Diophantine group over $\Q$.
\end{lemma}
\begin{proof}
Minkowski's bound $\sqrt{13}/2<2$ gives class number one for $\cO_F=\Z[d]$.
The prime $2$ is inert and $d$ has norm $-3$.
The units $-1,d+1$ realize all four sign patterns, so Dirichlet's unit theorem makes them generators of the unit squareclasses; in particular, all totally positive units are squares.

For $y^2=x^3+Bx$, the two-isogeny descent map is the squareclass of $x$, with values $B$ at $(0,0)$ and $1$ at $O$ \cite[Ch.~X, \S4]{Sil09}. At a prime not dividing $B$, a positive or negative valuation of $x$ gives respectively $v(y^2)=v(x)$ or $3v(x)$, so $v(x)$ is even. Thus its classes are supported on the primes dividing $B$. We apply this to $E$, where $B=-d^2$, and to the two-isogenous curve $y^2=x^3+4d^2x$, where $B=4d^2$.
For $B=-d^2$, its image is $\langle-1,d+1,d\rangle$: these classes are represented by $(0,0),Q,(d,0)$ and are independent by signs and the valuation at $(d)$.
For $B=4d^2$, all non-zero $x$-coordinates are totally positive and squares at $(d)$.
Indeed, $v_{(d)}(x)=1$ would give $v_{(d)}(y^2)=3$, since $-1$ is a nonsquare in $\F_3$; otherwise $x^2+4d^2$ is a square by Hensel's lemma.
The local square condition makes $v_{(d)}(x)$ even. The support bound and the fact that all totally positive units are squares therefore leave only the global classes $1,2$; the latter is a nonsquare modulo $(d)$.
The two image orders are $8$ and $1$; their product is $2^{2+\operatorname{rank}E(F)}$, so the rank is one.

At $(1-d)$ the good reduction is $y^2=x^3-x$ over $\F_3$, with four points, all reductions of $E[2]$.
The reduction kernel over $\Q_3$ is torsion-free by the formal logarithm \cite[Ch.~IV, Theorem~6.4; Ch.~VII, \S2]{Sil09}.
Hence the torsion is exactly $E[2]$. Membership in $\Gamma$ is expressed by $P=[2]R$ for some $R\in E(F)$. In the chosen rational coordinates, the doubling law makes this condition Diophantine.
\end{proof}

Fix a generator $P_0$ of $\Gamma$.

\begin{lemma}\label{lem:density}
The pairs $(a(P),b(P))$ with $P\in\Gamma\setminus\{O\}$ are dense in $\R^2$.
For every $P\in\Gamma\setminus\{O\}$, one has $0<u(P)<1/3$. Given any finite set of rational primes, $P$ may be chosen so that $u(P)$ is arbitrarily small at all of them.
\end{lemma}
\begin{proof}
Set $\lambda=\sqrt3/d$. Since $\lambda^2d=-d^\sigma$, the map
$\varphi(x,y)=(\lambda^2x,\lambda^3y)$ identifies $E$ with its conjugate $E^\sigma$ over $F(\sqrt3)$.
Being doubles, $\varphi(P_0)$ and $P_0^\sigma$ lie on the identity component of $E^\sigma(\R)$. They are independent modulo torsion: conjugating an integral relation by $\sqrt3\mapsto-\sqrt3$ and subtracting makes a multiple of $\varphi(P_0)$ torsion, so its coefficient, and then the other coefficient, is zero.
By Kronecker's theorem, this independence makes their simultaneous multiples dense in the product of the two circles.
On the identity component of $E_D(\R)$, the $y$-coordinate takes every real value away from $O$: the function $x^3-D^2x$ increases from $0$ to infinity for $x\ge|D|$.
Since $y(\varphi P)=\lambda^3y(P)$, the two conjugate values of $y(nP_0)$ vary densely in $\R^2$.
Their sum and difference recover $a(nP_0),b(nP_0)$ by an invertible linear change of coordinates, proving the first assertion.

A non-zero point of $\Gamma$ is a non-torsion double, so at each real embedding its $x$-coordinate exceeds the absolute value of the corresponding conjugate of $d$. Since $|\Nm_{F/\Q}(d)|=3$, this gives $0<u(P)<1/3$.
At all places of $F$ above the specified rational primes, choose open subgroups sufficiently close to $O$. Choose a non-zero multiple $P$ of $P_0$ lying in these subgroups. Then $u(P)$ is as small as required at every specified prime.
\end{proof}

\section{Integer arithmetic on the elliptic curve}\label{sec:model}

We are now ready to show that $n\mapsto nP_0$ is a Diophantine model of $\Z$ (Proposition~\ref{prop:model}). The group law adds indices, and we can recover multiplication from squaring, so our main task is to express squaring by a Diophantine condition.

\begin{proposition}\label{prop:model}
The bijection $n\mapsto nP_0$ is a Diophantine model of $\Z$ over $\Q$: it carries the graphs of addition and multiplication to Diophantine subsets of $\Gamma^3$.
\end{proposition}

For $P_1=n_1P_0$, $P_2=n_2P_0$, write
\[
\Delta=n_1^2-n_2.
\]
We seek a Diophantine condition equivalent to $\Delta=0$. It will be \eqref{eq:five-shifts}, which requires the condition $\mathcal C(P_1,P_2;\xi)$ of Definition~\ref{def:globalconditions} at five values of $\xi$.
Throughout this section $e=24$ is the exponent of Proposition~\ref{prop:height}, $n_0=1$ is the index of the reference point $P_0$, and $L=F(i)$, with $\sigma$ extended by $\sigma(i)=i$.
Since $[i](x,y)=(-x,iy)$, membership in $\Gamma+[i]\Gamma=\Z[i]P_0$ is expressed by $U=R_1+[i]R_2$ with $R_1,R_2\in\Gamma$. This group is therefore Diophantine. It is torsion-free, since a non-zero endomorphism has finite kernel, and its non-identity points have non-zero $x$- and $y$-coordinates.
Choose $U_j$ in this group and set $V_j=P_j-U_j$ for $0\le j\le2$.
Writing $U_j=[\gamma_j]P_0$ and $V_j=[\delta_j]P_0$ gives
\begin{equation}\label{eq:gaussian-indices}
\gamma_j+\delta_j=n_j.
\end{equation}

\subsection{Local congruences}\label{sec:local}
We begin with conditions at a single place $\p$ of $L$ above a prime $p$, and show that they imply that $p^e$ divides $\Delta$. Put $z=-x/y$ and $\rho=1/x$, both zero at $O$. Positive valuation of $\rho$ detects reduction to $O$ at the primes of good reduction considered below.
The formal logarithm of the invariant differential $dz/\sqrt{1+4d^2z^4}$, truncated at degree $e$, is
\begin{equation}\label{eq:truncated-log}
f(Z)=\sum_{r=0}^{5}\frac{(-1)^r}{4r+1}\binom{2r}{r}d^{2r}Z^{4r+1}.
\end{equation}
For $j=1,2$, define
\[
B_j=\frac{f(z(U_j))}{f(z(U_0))},\qquad
C_j=\left(\frac{f(z(V_j))}{f(z(V_0))}\right)^\sigma,
\]
with non-zero denominators.
For $p>25$ and a place $\p\mid p$ of $L$, normalize $v=v_\p$ by $v(p)=1$.
For $a,b\in L$, the congruence $a\equiv b\pmod{\p^e}$ means $v(a-b)\ge e$; neither term is initially required to be integral.
These primes are unramified in $L$ and are primes of good reduction for $E$ and $E^\sigma$.
Consider the conditions
\begin{equation}\label{eq:local-test}
\begin{gathered}
v(\rho(U_j)),\ v(\rho(V_j^\sigma))\ge1\quad(0\le j\le2),\\
B_1\equiv C_1,\qquad B_2\equiv C_2,\qquad B_1^2\equiv B_2
\pmod{\p^e}.
\end{gathered}
\end{equation}

\begin{lemma}\label{lem:local}
Let $p>25$ and let $\p\mid p$ be a place of $L$. Conditions \eqref{eq:local-test} imply $p^e\mid\Delta$.
If $\Delta=0$ and $p\equiv17\pmod{312}$, there exists a choice of auxiliary points satisfying \eqref{eq:local-test} at $\p$.
\end{lemma}
\begin{proof}
The first line puts the points in the formal groups of the good reductions.
There the logarithm preserves the valuation of $z$ and carries $[i]$ to multiplication by $i$.
For a non-zero formal-group point $Q$, truncation in \eqref{eq:truncated-log} satisfies
\[
\frac{f(z(Q))}{\log_E(z(Q))}\in1+\p^e.
\]
Indeed, $v(z(Q))\ge1$, the coefficient of $Z^m$ has denominator dividing $m$ locally, and $m-1-v_p(m)\ge e$ for all omitted degrees $m\ge25$.
For $j=1,2$, applying logarithms to $[\gamma_0]U_j=[\gamma_j]U_0$ and to the conjugate of $[\delta_0]V_j=[\delta_j]V_0$ gives
\[
\gamma_0B_j\equiv\gamma_j,\qquad \delta_0C_j\equiv\delta_j
\pmod{\p^e}.
\]
The relative errors multiply integral Gaussian indices. If $U_j=O$, then $\gamma_j=B_j=0$ and the congruence is exact; the case $V_j=O$ is similar.
Since $\gamma_0+\delta_0=n_0=1$, \eqref{eq:gaussian-indices} and $B_j\equiv C_j$ imply
\[
B_j=(\gamma_0+\delta_0)B_j\equiv\gamma_j+\delta_j=n_j
\pmod{\p^e}.
\]
Thus $B_j$ is integral, and the last congruence gives $p^e\mid\Delta$.

Now let $p>25$ and $p\equiv17\pmod{312}$, so $p\equiv1\pmod8$, $p\equiv2\pmod3$, and $p\equiv4\pmod{13}$.
Thus $p$ splits completely in $L$, and $d/d^\sigma$ is a nonsquare modulo $p$, since $dd^\sigma=-3$ and $\leg{-3}{p}=-1$.
The reduction has $j$-invariant $1728$ and is ordinary since $p\equiv1\pmod4$. Its endomorphism ring is an imaginary quadratic order containing $\Z[i]$, hence equals $\Z[i]$ \cite[Ch.~V, Example~4.5 and Theorem~3.1(b)]{Sil09}.
At the chosen place, its Frobenius is $[\pi]$, with $\pi=\pi_\p\in\Z[i]$ and $\pi\overline\pi=p$; its action on the rational two-torsion gives $\pi\equiv1\pmod{2\Z[i]}$.
Since $E^\sigma$ is the quadratic twist of $E$ by the nonsquare $d^\sigma/d$ at $\p$, its Frobenius is $[-\pi]$, with the Gaussian actions identified by the twisting isomorphism.
Take
\begin{equation}\label{eq:gaussian-choice}
U_j=\left[\frac{1-\pi}{2}\right]P_j,\qquad
V_j=\left[\frac{1+\pi}{2}\right]P_j,
\end{equation}
so that $U_j+V_j=P_j$. Write $P_j=[2]Q_j$ with $Q_j\in E(F)$. Then $U_j=[1-\pi]Q_j$ and $V_j^\sigma=[1+\pi]Q_j^\sigma$ lie in the formal groups, since $[1-\pi]$ and $[1+\pi]$ vanish on the reductions of $F$-points of $E$ and $E^\sigma$, respectively.
The multipliers are non-zero, so the denominator logarithms do not vanish; the relative truncation estimate gives the same conclusion for the denominators defining $B_j,C_j$.
The exact logarithm ratios are $n_j$; truncation proves \eqref{eq:local-test} when $\Delta=0$.
\end{proof}

\subsection{Combining the local tests}\label{sec:combine}
We now define a single Diophantine condition which guarantees that, for every prime $p$ dividing $M(\xi)$, the conditions \eqref{eq:local-test} hold at some place above $p$. Primes do not occur in its statement; they appear only in its verification.
For $\xi\in\Q$ with $G(\xi)\ne0$, set
\[
R=14G(\xi),\qquad H_R(X)=\frac{R}{R+X^2}.
\]
The factor $14$ is chosen so that $2R\in7G(\xi)\Q^{\times2}$; we deal with the primes $2$ and $7$ in the proof of Lemma~\ref{lem:model-completeness}.

\begin{definition}\label{def:globalconditions}
For $P_1,P_2\in\Gamma$ and $\xi\in\Q$, write $\mathcal C(P_1,P_2;\xi)$ if $G(\xi)\ne0$ and there exist six triples $(U_0,U_1,U_2)$ of points of $\Gamma+[i]\Gamma$, indexed by $1\le k\le6$, and $c_\xi,s_\xi,t_\xi,z_\xi\in\Q$ with the following properties. For the $k$-th triple, let $V_j=P_j-U_j$, let $B_j,C_j$ be the ratios of \S\ref{sec:local}, and put
\[
W_k=\prod_{j=0}^{2}H_R\bigl(\rho(U_j)\bigr)H_R\bigl(\rho(V_j^\sigma)\bigr)\cdot
\prod_{j=1}^{2}H_R\!\left(\frac{B_j-C_j}{R^{e-1}}\right)\cdot
H_R\!\left(\frac{B_1^2-B_2}{R^{e-1}}\right).
\]
With $w_0$ defined by the first equality below, require
\begin{equation}\label{eq:norm-test}
w_0=1-\Nm_{L/\Q}\left(\prod_{k=1}^6(1-W_k)\right),\qquad
c_\xi(R+w_0^2)^2=2R,
\end{equation}
and
\[
(1-c_\xi)(s_\xi^4-c_\xi^2)(t_\xi^4-c_\xi^2)=z_\xi^2.
\]
All quotients in the construction are required to be defined.
\end{definition}

\begin{lemma}\label{lem:globaldivisibility}
If $\mathcal C(P_1,P_2;\xi)$ holds, then
\begin{equation}\label{eq:combined-divisibility}
M(\xi)^e\mid\Delta.
\end{equation}
\end{lemma}
\begin{proof}
In an integer-valued valuation with $w=v(R)$ odd, the denominator terms of $H_R(X)$ have different valuations. Hence
\begin{equation}\label{eq:H-valuation}
v(H_R(X))=\max\{0,w-2v(X)\}.
\end{equation}
Let $p\mid M(\xi)$. Then $w=v_p(R)=v_p(G(\xi))$ is positive and odd.
All $H_R$-values are integral above $p$ by \eqref{eq:H-valuation}.
A first-line condition in \eqref{eq:local-test} fails only if the corresponding $\rho$-value has valuation at most zero.
A failed congruence has normalized difference of valuation at most $(e-1)(1-w)\le0$.
Either failure makes $W_k\equiv0\pmod{\p^w}$.
If all triples failed at all places above $p$, the product in \eqref{eq:norm-test}, and then its norm, would be $1$ modulo $p^w$.
Thus $v_p(w_0)\ge w$ and $v_p(c_\xi)=-w$, contrary to Lemma~\ref{lem:parity}.
Some triple therefore satisfies \eqref{eq:local-test} above $p$, proving \eqref{eq:combined-divisibility}.
The successful triple may depend on the prime and the place; six triples suffice for the six primes arising at each shift in Lemma~\ref{lem:model-completeness}.
\end{proof}

\subsection{The tests force equality}\label{sec:model-heights}
We now use heights to deduce $\Delta=0$ from the divisibility \eqref{eq:combined-divisibility}. This divisibility concerns the original indices, so we compare it with the heights of fixed large multiples of the points: passing to a multiple increases the canonical height by the square of the multiplier and leaves $\Delta$ unchanged.

For a positive integer $m_0$ and $j=1,2$, write
\[
x(m_0P_j)=x_{j,0}+x_{j,1}\sqrt{13},
\]
setting $x_{j,0}=x_{j,1}=0$ when $P_j=O$.
\begin{lemma}\label{lem:model-soundness}
There exists an integer $m_0\ge1$, independent of $P_1,P_2$, such that the existence of $s\in\Q$ satisfying
\begin{equation}\label{eq:five-shifts}
\mathcal C(P_1,P_2;s+\beta)\qquad
\bigl(\beta\in\{0,x_{1,0},x_{1,1},x_{2,0},x_{2,1}\}\bigr)
\end{equation}
implies $\Delta=0$.
\end{lemma}
\begin{proof}
With the convention $h(x(O))=0$, the canonical height satisfies, uniformly in $m_0,n_j$,
\begin{equation}\label{eq:amplified-height}
h(x(m_0P_j))=2m_0^2n_j^2\hh(P_0)+O_E(1)
\end{equation}
\cite[Ch.~VIII, \S9]{Sil09}.
If $\Delta\ne0$, Proposition~\ref{prop:height} and Lemma~\ref{lem:globaldivisibility} imply
\[
h(s+\beta)\le H M(s+\beta)^{24}\le H|\Delta|.
\]
The shift $\beta=0$ bounds $h(s)$; subtraction bounds $h(x_{j,0}),h(x_{j,1})$ and hence $h(x(m_0P_j))$ by $C|\Delta|$, with $C$ independent of $m_0$.
Summing \eqref{eq:amplified-height} gives
\[
m_0^2(n_1^2+n_2^2)\hh(P_0)\ll|\Delta|\le n_1^2+n_2^2.
\]
Since $\hh(P_0)>0$, this is impossible for a fixed sufficiently large $m_0$.
\end{proof}

Fix such an $m_0$ for the rest of the section. Each condition in \eqref{eq:five-shifts} has its own auxiliary variables.

\subsection{Equality supplies the tests}\label{sec:progressions}
Conversely, when $\Delta=0$ we show that \eqref{eq:five-shifts} can be satisfied; we use the Green--Tao theorem to find the primes needed for the local conditions.
\begin{lemma}\label{lem:model-completeness}
If $\Delta=0$, there exists $s\in\Q$ satisfying \eqref{eq:five-shifts}.
\end{lemma}
\begin{proof}
We choose the common shift so that all primes of positive odd valuation of $R$ are congruent to $17\pmod{312}$. The rational roots of $G$ allow these choices at all five shifts to come from one progression.
Choose a positive multiple $D_0$ of $14$ clearing the denominators of $x_{1,0},x_{1,1},x_{2,0},x_{2,1}$.
Choose $m\in\Z$ so that
$r_{\beta,\theta}=m+D_0(\beta-\theta)>2D_0$ for all five shifts $\beta$ and all six roots $\theta$ of $G$.
The primes congruent to $17$ modulo $312$ and exceeding $\max(D_0,25)$ have positive relative density.
By the Green--Tao theorem \cite[Theorem~1.2]{GT08}, there is a progression of length $1+\max_{\beta,\theta}r_{\beta,\theta}$, say $q_0,q_0+v_0,\ldots$, with $v_0>0$.
Put
\[
s=\frac{q_0+mv_0}{D_0v_0},\qquad
D_0v_0(s+\beta-\theta)=q_0+r_{\beta,\theta}v_0.
\]
For each $\beta$, the six values on the right are distinct primes exceeding both $D_0$ and $v_0$, hence prime to $D_0v_0$.
Moreover $s+\beta-\theta>2$, so at each shift $\xi=s+\beta$,
\[
R=\frac{14p_1\cdots p_6}{(D_0v_0)^6}>2.
\]
The only primes at which $R$ has odd valuation are $2,7,p_1,\ldots,p_6$; at the latter six its valuation is $1$.
In particular, $R+X^2\ne0$ for all $X\in L$, since $R$ has odd valuation at the split prime $p_1$.

For $1\le k\le6$, we choose the $k$-th triple so that it satisfies \eqref{eq:local-test} at a place above $p_k$; this is possible by Lemma~\ref{lem:local}.
There $v(R)=1$, and the arguments of all nine $H_R$-factors have positive valuation; hence each factor is $1$ modulo that place.
Thus $1-W_k$ has positive valuation there; all factors are integral at all places above $p_k$.
The norm in \eqref{eq:norm-test} is therefore divisible by $p_k$, so $w_0\equiv1\pmod{p_k}$ and $v_{p_k}(c_\xi)=1$.

We have $0<c_\xi\le2/R<1$.
At $\ell=2,7$, write $w=v_\ell(R)=1-6v_\ell(D_0v_0)\le-5$.
Since $w$ is odd, $v_\ell(R+w_0^2)=\min\{w,2v_\ell(w_0)\}\le w$. Hence
\[
v_\ell(c_\xi)=v_\ell(2)+w-2v_\ell(R+w_0^2)\ge v_\ell(2)-w>0,\qquad v_2(c_\xi)\ge6.
\]
Finally $c_\xi\in2R\Q^{\times2}$, so
$\sqf(c_\xi)=7p_1\cdots p_6\equiv7\pmod8$.
The valuation of $c_\xi$ is even at all primes outside $\{2,7,p_1,\ldots,p_6\}$. Thus $c_\xi$ has a square denominator.
Lemma~\ref{lem:parity} gives the required rational solution at each shift, proving \eqref{eq:five-shifts}.
\end{proof}

\begin{proof}[Proof of Proposition~\ref{prop:model}]
By Lemmas~\ref{lem:model-soundness} and~\ref{lem:model-completeness}, \eqref{eq:five-shifts} is equivalent to $n_1^2=n_2$. Its conditions are Diophantine: membership in $\Gamma+[i]\Gamma$ and multiplication by $m_0$ are Diophantine, and the remaining conditions are rational equations over $L$ with non-vanishing conditions. The coordinates of \S\ref{sec:curve} include the identity. Thus squaring is Diophantine.

Addition is the group law. If $S_1,S_2,S_3$ are obtained by squaring $mP_0,nP_0,(m+n)P_0$, then $P''+P''=S_3-S_1-S_2$ has the unique solution $P''=mnP_0$ in $\Gamma$, since $\Gamma$ has no non-zero $2$-torsion.
\end{proof}

\section{The integers among rational numbers}\label{sec:integers}

We can now pick out the integers among the rational numbers: we prove Theorem~\ref{thm:Z} and then deduce Corollary~\ref{cor:main}.

\begin{proof}[Proof of Theorem~\ref{thm:Z}]
For each fixed integer $C_0>1$, the conditions of Definition~\ref{def:intervals} are decidable on the non-zero indices of $P,P'$. Their coordinates can be computed exactly; the interval contains an integer precisely when $a\le\lfloor b\rfloor$, and the remaining inequalities are between rational numbers. The Davis--Putnam--Robinson--Matiyasevich theorem \cite{DPR61,Mat70} makes this relation on the indices Diophantine. Proposition~\ref{prop:model} transfers its definition to $\Gamma$, and projection to $(u,a,b)$ proves that $T$ is Diophantine.

Now fix $C_0$ with $\log C_0>H+\log16$, where $H$ is the constant in Proposition~\ref{prop:height}. Suppose \eqref{eq:final-definition} holds, and write $q=m/n$ in lowest terms with $n>0$. Since $u>0$, the second condition in \eqref{eq:final-definition} implies $g\ne0$. For $c=u^2/g$, the substitution
\[
(s_0,t_0,z_0)=(us,ut,u^4z/g^2)
\]
transforms that equation into the equation of Lemma~\ref{lem:parity}. Put $x=16q+8$. If $p\mid M(x)$, then $v_p(c)=2v_p(u)-v_p(g)$ is odd and therefore positive. Hence $M(x)\mid \Num(u)$. Since the denominator of $x$ is $n/\gcd(n,16)$, Proposition~\ref{prop:height} implies
\begin{equation}\label{eq:denominator-bound}
\log(n/16)\le h(x)\le H\Num(u)^{24}.
\end{equation}
Let $k\in[a,b]\cap\Z$. If $q\ne k$, then
\[
\frac1n\le |q-k|\le b-a<C_0^{-\Num(u)^{24}}<\frac1n.
\]
The last inequality follows from \eqref{eq:denominator-bound}, the choice of $C_0$, and $\Num(u)\ge1$. Hence $q=k$.

Conversely, suppose $q\in\Z$. Then $x=16q+8$ is divisible by $8$ and $|x|\ge8$, so $g>0$ and $g\equiv7\pmod8$. By Lemma~\ref{lem:density}, choose $P\in\Gamma\setminus\{O\}$ such that, for $u=u(P)$,
\[
v_2(u)\ge2,\qquad v_p(u)\ge v_p(g)\quad(p\mid g).
\]
The bound $0<u<1/3$ gives $0<c=u^2/g<1$. Moreover, $v_2(c)\ge3$, and $c$ has a square denominator. Its squarefree part is congruent to $7$ modulo $8$. Lemma~\ref{lem:parity}, followed by the inverse substitution, gives $s,t,z$ in \eqref{eq:final-definition}. By the real density in Lemma~\ref{lem:density}, there exists $P'$ with
\[
a(P')<q<b(P'),\qquad b(P')-a(P')<C_0^{-\Num(u)^{24}}.
\]
The integer $q$ itself belongs to this interval, so the resulting triple lies in $T$.
The Diophantine expression of rational inequalities completes the proof.
\end{proof}

\begin{proof}[Proof of Corollary~\ref{cor:main}]
A recursively enumerable subset of $\Q^r$ has a recursively enumerable set of numerator-denominator representatives in $\Z^{2r}$. Apply the theorem of Davis, Putnam, Robinson and Matiyasevich, followed by Theorem~\ref{thm:Z}, to these representatives. This proves (1); the converse follows by enumerating rational solutions.

Introducing variables for additions and multiplications expresses any polynomial system as quadratic equations. Their sum of squares is one quartic, proving (2). Apply the same construction to a fixed recursively enumerable non-recursive subset of $\Z$, keeping its integer argument as a parameter, to obtain (3).

For (4), the set of coefficient tuples of monic polynomials of fixed degree over $\Q$ having no rational root is recursive, hence defined by an existential formula by (1). The equivalence of this formula with the absence of a root is a first-order statement, so it holds in any field elementarily equivalent to $\Q$. An existential formula is preserved under field inclusions, so a polynomial without a root in $K$ cannot acquire one in $L$. Apply this to a minimal polynomial (see \cite[Proposition~23(d)]{Koe16}). Finally, Poonen's reduction \cite[Theorem~1.1(i)]{Poo09b} shows that a decision algorithm in (5) would decide rational solvability for arbitrary varieties, contradicting (3).
\end{proof}

\begin{remark}\label{rem:variables}
Fix $C_0$ sufficiently large for Theorem~\ref{thm:Z} and a quadratic presentation of $T$ with $r$ auxiliary variables. After the substitution $(s_0,t_0,z_0)=(us,ut,u^4z/g^2)$, \eqref{eq:final-definition} has a quadratic presentation using nine intermediate variables: three to compute $g$, one for $c$ with $cg=u^2$, two for $\sigma_1=s_0^2$ and $\sigma_2=t_0^2$, and three for $\tau_1=\sigma_1^2-c^2$, $\tau_2=\sigma_2^2-c^2$ and $\tau=\tau_1\tau_2$. The remaining equation is $(1-c)\tau=z_0^2$; here $g\ne0$ follows from $cg=u^2$ and $u>0$. Four further variables express $(q-a)(b-q)\ge0$ as a sum of squares, which is equivalent to $a\le q\le b$ since $a\le b$. Together with $u,a,b,s_0,t_0,z_0$, these give at most $r+19$ variables besides $q$; taking the sum of squares gives a quartic definition of $\Z$. We do not evaluate $r$.
\end{remark}

\end{document}